\newcommand{\rem}[1]{}
\documentclass{amsart}
\usepackage{amsfonts,amssymb,amsmath,amsthm,mathrsfs}
\newtheorem{thrm}{Theorem}[section]

\newtheorem{prop}[thrm]{Proposition}

\newtheorem{remark}[thrm]{Remark}

\newtheorem{definition}[thrm]{Definition}

\begin{document}

\author[C.~A.~Mantica, L.~G.~Molinari and U.~C.~De]{Carlo~Alberto~Mantica, Luca~Guido~Molinari and Uday~Chand~De}
\address{C.~A.~Mantica: Physics Department,
Universit\`a degli Studi di Milano,
Via Celoria 16, 20133 Milano, Italy and 
I.I.S. Lagrange, Via L. Modignani 65, 20161, Milano, Italy 
-- L.~G.~Molinari (corresponding author):
Physics Department,
Universit\`a degli Studi di Milano and I.N.F.N. sez. Milano,
Via Celoria 16, 20133 Milano, Italy,
-- U.~C.~De: Department of Pure Mathematics, University of Calcutta,
35 Ballygaunge Circular Road, Kolkata 700019, West Bengala, India}
\email{carloalberto.mantica@libero.it, luca.molinari@mi.infn.it, uc\underbar{\;}de@yahoo.com}
\title[GRW perfect-fluid space-times]{A condition for a perfect-fluid space-time\\ 
to be a generalized Robertson-Walker Space-Time}

\begin{abstract}
A perfect-fluid space-time of dimension $n\ge 4$ with 1) irrotational velocity vector field, 2) null
divergence of the Weyl tensor, is a generalised Robertson-Walker space-time with Einstein fiber. 
Condition 1) is verified whenever pressure and energy density are related
by an equation of state. The contraction of the Weyl tensor with the velocity vector field is zero.
Conversely, a generalized Robertson-Walker space-time with null divergence of the Weyl tensor is
a perfect-fluid space-time.
\end{abstract}
\date{28 september 2015}
\maketitle
\section{Introduction}
Standard cosmology is modelled on Robertson-Walker metrics for the high symmetry imposed on space-time by the cosmological principle
(spatial homogeneity and isotropy). A wide generalization are the "generalized Robertson-Walker spacetimes",
introduced in 1995 by Al\'{\i}as, Romero and S\'anchez \cite{Alias1,Ali95}:

\begin{definition}
An  $n$-dimensional Lorentzian manifold is a generalized Robertson-Walker space-time (GRW)
if locally the metric may take the form:
\begin{align}
ds^2 = -dt^2 + q(t)^2 \, g^*_{\alpha\beta}(x_2,\dots,x_n)dx^\alpha dx^\beta, \quad \alpha,\beta=2\dots n \label{metric}
\end{align}
that is, it is the warped product $(-1)\times q^2 \mathscr M^*$, where $\mathscr M^*$ is a $(n-1)$-dimensional Riemannian manifold. 
If $\mathscr M^*$ has dimension 3 and has constant curvature, the space-time is a Robertson-Walker space-time.
\end{definition}
\noindent
Such spaces include the Einstein-de Sitter space-time, the Friedmann cosmological models, the static Einstein 
space-time and the de Sitter space-time. They are the stage for treatment of small perturbations of the 
Robertson - Walker metric. 
We refer to the works by Romero et al. \cite{Romero}, S\'anchez \cite{San98,San99}, Gutierrez and Olea \cite{Gutierrez} 
for a comprehensive presentation of geometric properties and physical motivations.

Recently Bang-Yen Chen proved the following deep result \cite{Ban14}:
A Lorentzian manifold of dimension $n\ge 4$ is a GRW
 space-time if and only if it admits a time-like vector, $X^jX_j<0$, such that 
\begin{align}
\nabla_k X_j = \rho g_{kj}. \label{Chen}
\end{align}
According to Yano \cite{Yano1944}, a vector field $v$ is {\em torse-forming} if $\nabla_kv_j = \omega_k v_j + f g_{jk}$, 
where $f$ is a scalar function and $\omega_k$ is a 1-form. Its properties in pseudo-Riemannian manifolds 
were studied by Mike\u{s} and Rach\r{u}nek \cite{Mikes00, Mikes05}.
The vector is named {\em concircular} if $\omega_k$ is a gradient (or locally a gradient); in this case
$v$ can be rescaled to a vector $X$ with the property \eqref{Chen} \cite{Mikes00}.

Mantica et al. \cite{Mantica2015} proved two sufficient conditions for a Lorentzian manifold of 
dimension $n\ge 4$ to be a GRW space-time: the first one is the existence of a 
concircular vector such that $u^iu_i=-1$. 
The other sufficient condition restricts the Weyl and Ricci tensors: $\nabla_m C_{jkl}{}^m=0$ and $R_{ij}=B u_i u_j$
where $B $ is a scalar field and $u$ is a time-like vector field. \\

Lorentzian manifolds with a Ricci tensor of the form 
\begin{align}
R_{ij}=A g_{ij} + B u_iu_j, \label{ricci}
\end{align}
where $A $ and $B $ are scalar fields and $u_iu^i=-1$, are often named {\em
perfect fluid space-times}. 
It is well known that any Robertson-Walker space-time is a perfect fluid space-time \cite{ONeill}, and for
$n=4$ a GRW space-time is a perfect fluid if and only if it is a Robertson-Walker space-time.\\
The form \eqref{ricci} of the Ricci tensor is implied by Einstein's equation if the energy-matter 
content of space-time is a perfect fluid with velocity vector field $u$. 
The scalars  $A $ and $B $ are 
linearly related to the pressure $p$ and the energy density $\mu $ measured in the locally comoving inertial frame. 
They are not independent because of the Bianchi identity $\nabla^m R_{im} = \frac{1}{2} \nabla_i R$, which translates
into
\begin{align}
\nabla^m(B u_ju_m) = \tfrac{1}{2}\nabla_j  [(n-2)A -B ]. \label{Bianchi}
\end{align}
Geometers identify the special form \eqref{ricci} of the Ricci tensor as the defining property of quasi-Einstein 
manifolds (with any metric signature). The Riemannian ones were introduced by Defever and Deszcz in 1991 \cite{Defever} 
(see also \cite{Deszcz1998} and Chaki et al. \cite{Chaki}). In \cite{Des90} %(Lemma 4.1, Theorem 4.1 and Corollary 4.1) 
Deszcz proved that a 
quasi-Einstein Riemannian manifold with null Weyl tensor and few other conditions, is a warped 
product $(+1)\times q^2 \mathscr M^*$, where
 $\mathscr M^*$ is a $(n-1)$-dimensional Riemannian manifold of constant curvature.\\
%\cite{Deszcz1998} in the
%study of quasi-umbilical hypersurfaces of pseudo-Euclidean spaces \cite{21,26}
%and by 
Pseudo-Riemannian quasi-Einstein spaces arose in the study of exact solutions of Einstein's equations.
Robertson-Walker space-times are quasi-Einstein (see \cite{Bee96}, \cite{63} and references therein). 
%For recent results on quasi-Einstein manifolds we refer to \cite{11,27}. 

Shepley and Taub studied a perfect-fluid space-time in dimension $n=4$,
with
equation of state $p=p(\mu)$ and the additional condition that the Weyl tensor has null divergence, 
$\nabla_m C_{jkl}{}^m = 0$. They proved the following: the space-time is conformally flat $C_{jklm} = 0$, 
the metric is Robertson-Walker, the flow is irrotational, shear-free and geodesic \cite{She67}. \\
A related result was obtained by Sharma \cite{Sha93} (corollary p.3584): 
if a perfect-fluid space-time in $n=4$ with  $\nabla_m C_{jkl}{}^m = 0$
admits a proper conformal Killing vector, i.e. $\nabla_i X_j + \nabla_j X_i = 2\rho g_{ij}$, then it is conformally flat 
($C_{ijkl}=0$). In the framework of Yang's gravitational theory, Guilfoyle and Nolan proved that
a $n=4$ perfect fluid space-time with $p+\mu\neq 0$ is a Yang pure space (i.e.
$\nabla_m C_{jkl}{}^m=0$ and $\nabla_k R=0$) if and only if it is a Robertson-Walker 
space-time \cite{Guilfoyle}.\\
Coley proved that any perfect fluid solution of Einstein's equations satisfying a 
barotropic equation of state $p=p (\mu)$ and $p+\mu\neq 0$,
which admits a proper conformal Killing vector parallel to the fluid 4-velocity, 
is locally a Friedmann-Robertson-Walker model \cite{Coley}. \\
De et al. \cite{14} showed that $n=4$ conformally flat almost pseudo Ricci-symmetric space-times,  i.e. 
$\nabla_k R_{ij} =(a_k+b_k) R_{ij} + a_j R_{ik} + a_j R_{jk}$, are Robertson-Walker space-times.

Riemannian spaces equipped with a torse-forming vector field were studied by Yano as early as 1944
\cite{Yano1944}; his results were extended to pseudo-Riemannian spaces by Sinyukov \cite{Sinyukov1979}.
They showed that the existence of such a vector implies the following local shape of the metric:
$ds^2=\pm (dx^1)^2 + F(x_1,\dots, x_n) d\tilde s^2$, where $d\tilde s^2$ is the metric of the submanifold
parametrized by $x_2,\dots, x_n$. If the vector field is concircular (then it is rescalable to 
$\nabla_k X_j =\rho g_{kj}$) then $F$ is a function of $x_1$ only.

De and Ghosh \cite{De2000}
showed that if $R_{ij}=A g_{ij} + B u_iu_j$ with $u_i$ closed and $C_{ijkl}=0$, then $u$ is a 
concircular vector. The results were extended by Mantica et al. to pseudo Z-symmetric spaces \cite{Mantica2012} and to
weakly Z-symmetric spaces \cite {Mantica2012b}.\\

In this paper the theorem by Shepley and Taub is generalised to perfect-fluid space-times of dimension $n\ge 4$. The converse is also proven: a GRW space-time with $\nabla_m C_{jkl}{}^m=0$ is a perfect-fluid space-time.
In the conclusion, some consequences for physics are presented.

\section{The theorem}
\begin{thrm}\label{THRM1}
Let $\mathscr M$ be perfect fluid-space-time, i.e. a Lorentzian manifold (of dimension $n>3$) 
with Ricci tensor  $R_{kl} =A g_{kl}+B u_k u_l $, 
where $A $ and $B $ are scalar fields, $u$ is a time-like unit vector field $u^j u_j =-1$.\\
If  $\nabla_k u_j -\nabla_j u_k=0$ ($u$ is closed) and if $\nabla_m C_{jkl}{}^m =0$, then:\\
i) $u$ is a concircular vector and it is rescalable to a time-like conformal 
Killing vector $X$ such that
\begin{align}
\nabla_k X_j = \rho g_{kj} \quad \text{and} \quad \nabla_k \rho = \frac{A-B}{1-n}X_k;
\end{align} 
ii) $\mathscr M$ is a generalised Robertson-Walker space-time whose sub-manifold $(\mathscr M^*, g^*)$ is 
a Riemannian Einstein space.\\
iii) $C_{jklm}u^m =0$.
\begin{proof}
The condition $\nabla_m C_{jkl}{}^m =0$ implies:
$\nabla_k R_{jl}-\nabla_l R_{jk} = \frac{1}{2(n-1)}(g_{jl} \nabla_k R - g_{jk}\nabla_l R)$.
With the explicit form of the Ricci tensor, it becomes
\begin{align}
\nabla_k (B u_j u_l) -\nabla_l (B u_j u_k) 
= - \frac{g_{jl} \nabla_k \gamma- g_{jk}\nabla_l \gamma }{2(n-1)}  \label{1.1}
\end{align}
being $\gamma = (n-2)A +B $. By transvecting with $u^ju^l$ and using 
$u^l\nabla_k u_l=0$ we obtain
\begin{align}
(\nabla_k + u_ku^l\nabla_l )B +B u^l\nabla_l u_k = \frac{1}{2(n-1)}(\nabla_k +
u_ku^l\nabla_l )\gamma . \label{1.3}
\end{align}
Contraction of the identity \eqref{Bianchi} with $u^j$ gives $-B \nabla_m u^m = \frac{1}{2}u^m \nabla_m \gamma $, which rewrites 
identity \eqref{Bianchi} as:
\begin{align}
(\nabla_k +u_ku^i\nabla_i)B +B u^m\nabla_m u_k  = \tfrac{1}{2}(\nabla_k +u_ku^i\nabla_i)\gamma. \label{Bianchi3}
\end{align}
Equations \eqref{1.3} and \eqref{Bianchi3} imply:
\begin{align}
(\nabla_j +u_j u^k\nabla_k)\gamma =0, \label{1.5} \\
(\nabla_j + u_ju^k\nabla_k ) B+ B u^m\nabla_m u_j =0. \label{1.6}
\end{align}
%Eq.\eqref{1.5}  implies
%\begin{align}
%u_i\nabla_j \gamma =  u_j \nabla_i \gamma  \label{1.7}
%\end{align}
Contraction of \eqref{1.1} with $u^l$ gives:
\begin{align*}
& -u_j (\nabla_k +u_ku^l\nabla_l)B - B\nabla_k u_j -u_j Bu^l\nabla_l u_k - u_kB u^l\nabla_l u_j \nonumber\\ 
& =-\frac{1}{2(n-1)} (u_j\nabla_k -g_{jk}u^l\nabla_l)\gamma 
\end{align*}
By use of eq. \eqref{1.6} it simplifies to:
\begin{align}
 B(\nabla_k  +u_k u^m\nabla_m )u_j  =\frac{1}{2(n-1)} (u_j\nabla_k -g_{jk}u^l\nabla_l)\gamma \label{betabeta}
\end{align}
If $u$ is closed it is $u^m\nabla_m u_j = u^m\nabla_j u_m =0$. Eq.\eqref{betabeta} simplifies and shows that
$u$ is a torse-forming vector:
\begin{align}
 \label{1.11}
\nabla_k u_j = \frac{\nabla_k \gamma }{2B (n-1)} u_j - \frac{u^m\nabla_m \gamma}{2B(n-1)} g_{kj} \equiv
 \omega_k u_j+ f g_{kj}
\end{align}
Let us show that $u$ is a concircular vector, i.e. that $\omega_k$ is closed:\\ 
$\nabla_j \omega_k-\nabla_k\omega_j =  -\frac{1}{B}  (\omega_k\nabla_j-\omega_j\nabla_k) B = -(\omega_k u_j -\omega_j u_k) u^m\nabla_m B$
by \eqref{1.6}.  Eq. \eqref{1.5} gives the relation $\omega_k = - u_k u^m \omega_m$, then $\omega_ku_j-\omega_j u_k=0$.

Being closed, $\omega_k$ is locally the gradient of a scalar function: $\omega_k=\nabla_k \sigma $. %(see Hall \cite{Hall} page 242-243). 
Let 
$X_l = u_l e^{-\sigma}$; we have $\nabla_k X_l = e^{-\sigma} (-u_l \nabla_k \sigma +\omega_k u_l+ f g_{kl}) = e^{-\sigma} f g_{kl} $ and consequently
\begin{align}
\nabla_k X_l =\rho \, g_{kl} \label{eqforchen}
\end{align}
being $\rho=e^{-\sigma}f$ and $X_j X^j=-e^{-2\sigma} < 0$  (time-like vector). The symmetrized equation $\nabla_k X_j +\nabla_j X_k = 2\rho g_{kj}$ shows that $X_j$ is a 
conformal Killing vector \cite{63}.

According to Chen's theorem, \eqref{eqforchen} is a sufficient condition for the space-time to be a GRW. In appropriate coordinates 
$\mathscr M = (-1)\times q^2 \mathscr M^*$. The additional condition $\nabla_m C_{jkl}{}^m=0$ assures that the $(n-1)$-dimensional 
Riemannian space $\mathscr M^*$ is an Einstein space, by G\c{e}barowski's lemma \cite{Gebarowski}.

Another derivative and the Ricci identity give: $(\nabla_j\nabla_k-\nabla_k\nabla_j) X_l =R_{jkl}{}^m X_m = 
g_{kl}\nabla_j\rho - g_{jl}\nabla_k \rho $. Contraction with $g^{kl}$: $R_{jm} X^m = (1-n) \nabla_j \rho $. However, for the perfect
fluid \eqref{ricci},
$R_{jm}X^m = (A -B)X_j$, then:
\begin{align}
\nabla_j \rho = \frac{A-B}{1-n} X_j \label{CKV}
\end{align}
(this is an explicit expression for a relation obtained by Chen). 
Therefore, if $A\neq B$ the conformal killing vector $X$ is proper; if $A=B$ it is homothetic. Moreover:
\begin{align}
R_{jklm} X^m = \frac{A-B}{1-n} (X_j g_{kl} - X_k g_{jl})
\end{align}
The Weyl tensor is:
\begin{align*}
C_{jklm}=R_{jklm} +\tfrac{1}{n-2}(g_{jm}R_{kl}-g_{km}R_{jl}+R_{jm}g_{kl}-R_{km}g_{jl}) - 
\frac{(g_{jm}g_{kl}-g_{mk}g_{jl})R}{(n-1)(n-2)}
\end{align*}
The previous equations and little algebra imply that $C_{jklm}X^m=0$, so that $C_{jkl}{}^mu_m=0$. It follows that 
the Weyl tensor is purely electric \cite{Ortaggio2013}. \\
In $n=4$ the condition is equivalent to $u_i C_{jklm}+ u_j C_{kilm}+ u_k C_{ijlm}=0 $ (see Lovelock and 
Rund \cite{LR} page 128). Multiplication by $u^i$
gives $C_{ijkl}=0$.
\end{proof}
\end{thrm}

\begin{remark}
Eq.\eqref{1.5} gives $u_j\nabla_k\gamma =u_k\nabla_j\gamma $. In the antisymmetric part of eq.\eqref{betabeta},  
$B(\nabla_k u_j - \nabla_j u_k) + B( u_ku^m\nabla_m u_j -  u_j u^m\nabla_m u_k) = 0$, the last terms are replaced 
with the help of \eqref{1.6} to give:
\begin{align}
\nabla_k (Bu_j) = \nabla_j (B u_k) \label{BU}
\end{align}
\end{remark}

\begin{remark}
The case $A =0$, i.e. $R_{ij}=B u_i u_j$,  was studied in \cite{Mantica2015}. Since $\gamma =-B$, the property 
$u_j\nabla_k\gamma =u_k\nabla_j\gamma $ and \eqref{BU} imply that $u$ is closed.\\
If $A \ne 0$ the condition that $u$ is closed is necessary for proving the theorem. However, if a
one-to-one differentiable relation $A (x)= F(B (x))$ exists, one proves that
$u$ is closed.
\end{remark}

\begin{remark}
In \cite{Broz12,Cat11,Fer12} a notion of quasi-Einstein manifold different from \eqref{ricci} was introduced. It emerges from 
generalizations of Ricci solitons. More generally, they defined a generalized quasi - Einstein manifold by the condition 
\begin{align}
R_{ij} +\nabla_i \nabla_j \theta -\eta (\nabla_i\theta)(\nabla_j\theta) = \lambda g_{ij}  \label{1.16}
\end{align}
where $\theta, \eta, \lambda$  are smooth functions. If  $\lambda=$ const and $\eta=0$ it is named gradient Ricci soliton, 
if $\lambda=$const. and $\eta =$const. it is named quasi-Einstein.\\ 
In the present case, the condition that $u$ is closed means that locally $u_k =
\nabla_k \theta$, for some function $\theta$. Then \eqref{ricci} takes the form 
$R_{ij} = A g_{ij} + B (\nabla_i\theta)(\nabla_j\theta)$. At the same time, eq.\eqref{1.11}
can be written 
$\nabla_i \nabla_j \theta = f(\nabla_i\theta)(\nabla_j\theta) + fg_{ij}$ (since $f=-u^k\omega_k$ and $\omega_i = - u_iu^k\omega_k$ by eq.\eqref{1.5},
it is $\omega_i = f u_i$). The sum of the equations yields  a Ricci tensor of the form  \eqref{1.16}
with $\lambda = A+f$ and  $\eta=B+f$, i.e. the manifold is generalized quasi-Einstein in the sense of \cite{Broz12,Cat11,Fer12}.
A gradient Ricci soliton is recovered if $A +f=$const. and $B+f=0$.\\ 
In \cite{Broz12} it was proven that locally conformally flat Lorentzian quasi-Einstein manifolds are globally 
conformally equivalent to a space form, or locally isometric to a warped product of Robertson-Walker type, or a pp-wave. 
Catino \cite{Cat11} proved that a complete (i.e. $A +f $ is a smooth function) 
generalized quasi-Einstein Riemannian manifold with harmonic Weyl tensor and 
zero radial curvature, is locally a warped product with $(n-1)$ dimensional Einstein fibers.  
\end{remark}

An inverse statement of the theorem is proven:

\begin{thrm}
A generalized Robertson-Walker space-time with $\nabla_m C_{jkl}{}^m = 0$ is a quasi-Einstein space-time. 
\begin{proof}
A GRW is characterized by the metric \eqref{metric}. The explicit form of the Ricci tensor $R_{ij}$ is reported for example in Arslan et al.\cite{Des2014}:
$R_{1\alpha}=R_{\alpha 1}=0$,
\begin{align*}
R_{11} = - (n-1) \frac{q'}{q},  \quad R_{\alpha\beta} = R^*_{\alpha\beta} + g^*_{\alpha\beta}\left[ q^{\prime 2} (n-2)+ qq^{\prime\prime}\right ], \quad
\alpha,\beta=2\dots n.
\end{align*}
G\c{e}barowski proved that $\nabla_m C_{jkl}{}^m = 0$ if and only if $R^*_{\alpha\beta} = g^*_{\alpha\beta} \frac{R^*}{n-1}$, then:
\begin{align*}
 R_{\alpha\beta} =  g^*_{\alpha\beta}\left[\frac{R^*}{n-1} + \, q^{\prime 2} (n-2) + qq^{\prime\prime}\right ].
\end{align*}
Following the trick in \cite{11}, in the local frame where \eqref{metric} holds, 
define the vector $u^1=1$ and $u^\alpha =0$ ($u_1=-1 $). It is $u_ju^j = -1 $ in any frame. The components of the Ricci tensor gain the covariant expression $R_{ij}= Ag_{ij}+ B u_iu_j$, where:
\begin{align}
A=\frac{1}{q^2}\left[\frac{R^*}{n-1} + \, q^{\prime 2} (n-2) + qq^{\prime\prime}\right ], \quad B=-(n-1)\frac{q^\prime}{q}+A
\end{align}
The expression is such in all coordinate frames, and characterizes a quasi-Einstein Lorentzian manifold.
\end{proof}
\end{thrm}

\section{Some notes on physics}
We transpose some of the results to physics (we use units $c=1$). Consider a perfect fluid with energy 
momentum tensor $T_{ij} = p g_{ij} + (p+\mu) u_iu_j$,
where $u_j$ is the velocity vector field, $p$ is the isotropic pressure field and $\mu $ is the energy density. By 
Einstein's equations $R_{ij}-\frac{1}{2} Rg_{ij} = \kappa T_{ij}$ ($\kappa = 8\pi G$ 
is the gravitational constant) the Ricci tensor is:
$$ R_{ij} = \kappa (p+\mu) u_iu_j + \kappa \frac{p-\mu}{2-n} g_{ij} . $$ 
Comparison with the form \eqref{ricci} identifies  $A = \kappa (p-\mu)/(2-n)$, $B = \kappa (p+\mu)$. Then $\gamma = 
(n-2)A +B = 2\kappa \mu$.\\
As is well known (see Wald \cite{Wald}) in General Relativity the equations of motion  
$\nabla_k T^{kj}=0$  result from the Bianchi identity in Einstein's equations. 
For a perfect fluid, the projection along $u$ and its complementary part are:
\begin{align}
& u^k\nabla_k \mu + (p+\mu) \nabla_k u^k =0 \label{eqm1}\\
& (\nabla_j  + u_ju^k\nabla_k) p + (p+\mu) u^k\nabla_k u_j = 0 \label{eqm2}
\end{align}
%We show that if a one-to-one constitutive relation $p=p(\mu)$ is given, $p+\mu\neq 0$, and $\nabla_m C_{jkl}{}^m=0$,
%then the integral lines of the motion are geodesics and the velocity vector field is irrotational (i.e. closed). 
By taking into account the results of the previous section we prove:

\begin{prop}
A perfect fluid space-time in dimension $n\ge 4$, with differentiable equation of state $p=p(\mu)$, $p+\mu \neq 0$, and with 
null divergence of the Weyl tensor, $\nabla_m C_{jkl}{}^m =0$, is a generalized Robertson-Walker space-time. \\
The velocity vector field is irrotational ($\nabla_k u_l-\nabla_l u_k=0$), geodesic ($u^k\nabla_k u_j=0$) and it annihilates the Weyl tensor 
($C_{jkl}{}^m u_m=0$). 
\begin{proof}
We prove that $u$ is irrotational and geodesic. Then, by the main theorem \ref{THRM1} it follows that the manifold is a generalized
Robertson-Walker space-time  and that $u$ annihilates the Weyl tensor.\\
If $p'(\mu)\neq 0$ then $\nabla_k p = p'(\mu) \nabla_k\mu $. The eqs. $u_j\nabla_k \gamma = u_k\nabla_j\gamma$ and \eqref{BU} become:
$u_j\nabla_k\mu = u_k\nabla_j\mu$ and $\nabla_j [(p+\mu)u_k] = \nabla_k [(p+\mu) u_j]$. Being $\nabla_k p = p'(\mu)\nabla_k\mu$
it follows that $\nabla_k u_j=\nabla_j u_k$.\\
Eq.\eqref{1.5} is $\nabla_j \mu + u_ju^m\nabla_m \mu=0$, and translates to $\nabla_j p + u_ju^m \nabla_m p=0$. This is used in \eqref{eqm2}
to annihilate the first term. The equation of a geodesic is obtained: $(p+\mu) u^k\nabla_k u_j = 0$. If $\nabla_k p=0$, 
eq. \eqref{eqm1} again gives $(p+\mu)u^k\nabla_k u_j=0$.\\
\end{proof}
\end{prop}

The special case $A=B$ in \eqref{CKV} characterizes a homothetic conformal Killing field ($ \nabla_j X_k=\rho g_{jk}$ with $\nabla_j\rho =0$).
In terms of pressure and density this means 
$$p=\frac{3-n}{n-1}\mu $$ 
which, in $n=4$, is $p=-\mu /3$.

\end{document}